\documentclass[final, leqno]{amsart}

\usepackage{amsmath,amsthm} %,amsfonts}
\usepackage {latexsym}
\usepackage{amssymb}
\newcommand{\eps}{{\varepsilon}}
\newcommand{\R}{{\mathbb R}}

\newcommand{\Compl}{{\mathbb C}}

\newcommand{\B}{{\mathcal B}}
\newcommand{\W}{{\mathcal W}}

\newcommand{\les}{\lesssim}

\newcommand{\Kato}{{\mathcal K}}
\newcommand{\txt}{\textstyle}
\newcommand{\scr}{\scriptstyle}

\newcommand{\Lap}{\Delta}

\newcommand{\Res}{R_0^+}

\newcommand{\vp}{\varphi}

\newcommand{\la}{\langle}
\newcommand{\ra}{\rangle}

\def\norm[#1][#2]{\|#1\|_{#2}}
\def\bignorm[#1][#2]{\big\|#1\big\|_{#2}}
\def\Bignorm[#1][#2]{\Big\|#1\Big\|_{#2}}
\def\japanese[#1]{\langle #1 \rangle}
\def\Im[#1]{{\rm Im}(#1)}
\def\Re[#1]{{\rm Re}(#1)}

\newcommand{\be}{\begin{equation}}
\newcommand{\ee}{\end{equation}}

\newtheorem{theorem}{Theorem}

\theoremstyle{remark}
\newtheorem{remark}{Remark}

\begin{document}

\title[Dispersive Bound for Scaling-Critical Schr\"odinger Operators]{Schr\"odinger Dispersive Estimates for a Scaling-Critical Class of Potentials}
\date{September 12, 2010}

\author{Marius\ Beceanu}
\thanks{The first author was supported by the ANR program PREFERED in the process of writing this paper}
\address{EHESS--CAMS, 54 bd Raspail, Paris 75006, France}
\email{mbeceanu@ehess.fr}

\author{Michael\ Goldberg}
\thanks{The second author received support from NSF grants
DMS-0600925 and DMS-0901063 during the preparation of this work.}
\address{Department of Mathematics, University of Cincinnati,
Cincinnati, OH 45221-0025}
\email{Michael.Goldberg@uc.edu}

\begin{abstract}
We prove a dispersive estimate for the evolution of Schr\"odinger operators
$H = -\Delta + V(x)$ in ${\mathbb R}^3$.  The potential should belong
to the closure of $C^c_b(\R^3)$ with respsect to the global Kato norm.
Some additional spectral conditions are imposed, namely that no resonances
or eigenfunctions of $H$ exist anywhere within the interval $[0,\infty)$.
The proof is an application of a new version of Wiener's $L^1$-inversion
theorem. 
\end{abstract}

\maketitle

\section{Introduction} \label{sec:intro}

Solutions to the linear Schr\"odinger equation
are governed by a number of dispersive and smoothing estimates.    These
inequalities place limits on the types of singularities that can arise
as well as the length of time they are allowed to persist.  On short time
scales the dispersive bounds are a useful
stepping stone toward a nonlinear local existence theory, and on
long time scales they enable analysis of asymptotic properties and scattering
behavior.  We will examine the use of initial data in $L^1(\R^3)$ to control
the supremum norm of the solution at later times.  Mappings between those
spaces are described fully by the pointwise size of the propagator kernel
without regard to oscillations in sign.

To see that oscillatory integrals play a major role at every other step of the
computation, consider the case of the free Schr\"odinger equation.
Intial data at time zero is brought forward to time $t$ through the
action of a Fourier multiplier $e^{it\Delta}$. In spatial variables
this is equivalent to convolution against the complex Gaussian kernel 
$(-4\pi it)^{-3/2} e^{i(|x|^2/4t)}$.  The decay rate of $|t|^{-3/2}$
(more generally $|t|^{-n/2}$ when the equation is set in $\R^n$)
arises as a consequence of stationary phase within the Fourier inversion
integral. 
It immediately follows
that the free evolution satisfies a dispersive bound 
\begin{equation} \label{eq:freedispersive}
\norm[e^{it\Lap}f][\infty] \le (4\pi|t|)^{-3/2}\norm[f][1]
\end{equation}
 at all times $t \not= 0$.
In this paper we seek to prove similar estimates for the time evolution
$e^{-itH}$ induced by a perturbed Hamiltonian $H= -\Delta + V(x)$.
The class of admissible potentials will be defined with respect to the global
Kato norm
\begin{equation}
\norm[V][\Kato] = \sup_{y\in\R^3} \int_{\R^3} \frac{|V(x)|}{|x-y|}\,dx.
\end{equation}
In relation to other usual classes of functions, $L^{3/2-\epsilon} 
\cap L^{3/2+\epsilon} \subset L^{3/2, 1} \subset \Kato$, where $L^{3/2, 1}$
 is a Lorentz space.

We assume that $V$ belongs to $\Kato_0 \subset \Kato$, the norm-closure
of the bounded, compactly supported functions within $\Kato$.  This subspace
carries the Kato norm, which is homogeneous with respect to the scaling
$V_r(x) = r^2 V(rx)$, $r > 0$.  There are no further restrictions on the size of 
$V$ or
its negative or imaginary parts.

Dispersive estimates cannot hold for all initial data if $H$ possesses 
one or more bound states.  Whenever there exists a nonzero function 
$\Psi \in L^1(\R^3)$ that solves the eigenvalue equation 
$H\Psi = \lambda\Psi$, the associated Schr\"odinger evolution
\begin{equation*}
e^{-itH}\Psi = e^{-it\lambda}\Psi
\end{equation*}
maintains a constant amplitude for all times
in violation of~\eqref{eq:freedispersive}.  
Bound states are typically removed from consideration by applying a spectral
projection to the initial data.  A revised dispersive estimate for 
$H = -\Delta + V$ might take the form
\begin{equation} \label{eq:dispersive1}
\big\|e^{-itH}(I- {\txt \sum}_j P_{\lambda_j}(H))f\big\|_\infty 
  \leq |t|^{-\frac32}\norm[f][1]
\end{equation}
where $P_{\lambda_j}(H)$ is a projection onto the point spectrum of $H$
at the eigenvalue $\lambda_j$.  

One additional concern here is the possible existence of resonances, which
are solutions of the equation
$\psi + (-\Delta - \lambda \pm i0)^{-1} V \psi = 0$ that do not decay
rapidly enough to belong to $L^2(\R^3)$ but instead satisfy 
$(1+|x|)^{-s} \psi\in L^2$ for every $s > \frac12$.  
Applying $-\Delta -\lambda$ to both sides shows that $H\psi = \lambda\psi$
for these functions as well.
Resonances exhibit enough
persistence behavior (by virtue of their resemblance to $L^2$ bound states) to
negate most dispersive estimates, but they cannot be so easily removed with a 
spectral projection.  

\begin{theorem}
\label{thm:dispersive}
Let $V \in \Kato_0$ be a real-valued potential for which the
Schr\"odinger operator $H = -\Delta + V$ has no resonances or eigenvalues
within the halfline $[0, +\infty)$.  Then
\begin{equation} \label{eq:dispersive}
\big\|e^{-itH} (I - P_{pp}(H))f\big\|_{\infty} 
  \les |t|^{-\frac32}\norm[f][1]. 
\end{equation}
Under these conditions, the spectral multiplier $P_{pp}(H) := 
\sum_j P_{\lambda_j}(H)$ is a finite-rank projection.
\end{theorem}

\begin{remark}
No part of the proof of Theorem~\ref{thm:dispersive} relies on the potential
$V(x)$ being real-valued or on $H$ being self-adjoint (the use of spectral
measures is just a convenience).  As stated the
theorem is equally true for complex $V\in \Kato_0$ so long as $H$ satisfies the
assumed spectral conditions.  In the complex case the projections 
$P_{\lambda_j}(H)$ should be constructed to take into account the entire
generalized eigenspace over $\lambda_j$.  Even so they are still finite-rank
operators that can be recovered from the (analytic) functional calculus of $H$
in the neighborhood of each $\lambda_j \in \sigma(H)$.
\end{remark}

Based on the commutation relation between dilations and the Laplacian,
%$D_{r^{-1}} \Delta D_r = r^{-2} \Delta$,
a potential of the form
$V_r(x) = r^2 V(rx)$ is guaranteed to produce the same dispersive bounds
as $V$ itself.  For this reason we regard an inverse-square law for
potentials as being critical with respect to scaling.  In the explicit
case $V(x) = C|x|^{-2}$, dispersive bounds are true only when 
$C \ge 0$~\cite{BuPlStTa03}.  If the pointwise decay rate is further relaxed
to any lower power $|V(x)| \le C|x|^{-2+\eps}$, the dispersive bound may
fail even for nonnegative potentials~\cite{GoVeVi06}.  We note that
the global Kato norm, and membership in the class $\Kato_0$, are preserved
among all members of a family $V_r$, $r>0$.

The first dispersive estimate of the form~\eqref{eq:dispersive} was proved
in~\cite{JoSoSo91} for real potentials satisfying
a regularity hypothesis $\hat{V} \in L^1(\R^3)$ and pointwise decay bound
$|V(x)| \les \japanese[x]^{-7-\eps}$.  Successive improvements (\cite{Ya95},
~\cite{GoSc04a}, and~\cite{Go06b} in chronological order)
relaxed the requirements
on $V$ down to the condition $V \in L^{\frac32 - \eps}(\R^3) \cap
L^{\frac32 + \eps}(\R^3)$.  In terms of homogeneous functions, this permits
local singularities on the order of $|x|^{-2+\eps}$ and decay at the rate
$|x|^{-2-\eps}$ for large $x$.  In all these results $H$ is assumed not to
possess an eigenvalue or resonance at zero.

Each of the above conditions leads to a situation where
the measured "size" of a rescaled potential $V_r$ increases without
bound as $r \to 0$ and $r \to \infty$.  Prior to the current work the only
scale-invariant class of potentials known to produce a dispersive estimate is
based on a smallness condition $\norm[V][\Kato] < 4\pi$
that makes the perturbation series absolutely convergent~\cite{RoSc04}.
Eigenvalues and resonances are not possible (at zero or elsewhere) for
such small potentials.

The proof of Theorem~\ref{thm:dispersive} is based on a broad extension
of Wiener's $L^1$-inversion theorem to operator-valued functions,
first observed in~\cite{Be09p}.  In the one-dimensional setting, 
scaling-critical dispersive estimates were established with the help of
the classical (scalar) Wiener inversion theorem instead~\cite{GoSc04a}. 
There it is invoked at a crucial juncture to show that a
particular quotient of functions has integrable Fourier transform.  Keeping
the denominator nonzero ends up being equivalent to the absence of
resonances.

Our generalized inversion theorem plays a similar role here with operator
inverses taking the place of a quotient.  The precise statement and proof,
which contain the result in~\cite{Be09p} as a special case, are given in
Section~\ref{sec:Wiener}.  Once again the spectral property
required to apply the theorem to Schr\"odinger's equation is contingent
on keeping the continuous spectrum of $H$ free from eigenvalues and
resonances.

The mathematical argument divides neatly into three parts.
Section~\ref{sec:reduction} outlines the decomposition of
$e^{-itH}$ according to the spectral measure of $H$ and reduces the
dispersive bound to a desired estimate in $L^1$.  Section~\ref{sec:Wiener}
introduces the machinery related to Wiener's inversion theorem, and
in the concluding section we verify that the desired estimate fits into
this abstract framework.  

The same basic method also applies to the wave equation propagators
$\cos(t\sqrt{H})$ and  $\frac{\sin(t\sqrt{H})}{\sqrt{H}}$, yielding estimates 
in a variety of $L^p$ and Sobolev spaces.
The details of these cases will be presented in a separate paper~\cite{BeGo10p}
so that the present exposition can focus on operator-theoretic concerns
with a relatively small number of side calculations.

\section{Reduction to Resolvent Estimates} \label{sec:reduction}

The derivation of dispersive estimates from properties
of the free and perturbed resolvents is a standard practice
for time-independent Schr\"odinger operators on $\R^n$.  In fact it
is formally equivalent to methods based on the Duhamel propagation
formula, under the Fourier duality pairing of the time variable $t$ with
the spectral parameter $\lambda$.  The most prominent feature in
Duhamel's formula is a convolution integral (in $t$) against the free
propagator kernel.  In the dual setting this appears instead as a pointwise
operator acting on functions of $\lambda$.

Let $H=-\Delta+V$ in $\R^3$ and for each $z \in \Compl \setminus \R^+$
define the resolvents
$R_0(z) :=(-\Delta-z)^{-1}$ and $R_V(z):= (H-z)^{-1}$.
The operators $R_0(z)$ are all bounded on $L^2(\R^3)$ and
act explicitly by convolution with the kernel
\begin{equation*}
R_0(z)(x,y) = \frac{e^{i\sqrt{z}|x-y|}}{4\pi|x-y|},
\end{equation*}
where $\sqrt{z}$ is taken to have positive imaginary part.
While $R_V(z)$ is not translation--invariant and
there is no simple formula for its integral representation,
it can be expressed in terms of $R_0(z)$ via the identity
\begin{equation} \label{eq:ResIdent}
R_V(z) = (I + R_0(z)V)^{-1}R_0(z) = R_0(z)(I + VR_0(z))^{-1}. 
\end{equation}
In the case where $z = \lambda \in \R^+$, the resolvent may be defined
as a limit of the form $R_0(\lambda \pm i0) := \lim_{\eps\downarrow 0} 
R_0(\lambda\pm i\eps)$.  The choice of sign determines which branch of the
square-root function is selected in the formula above, therefore the two
continuations do not agree with one another.  
In this paper we refer to resolvents along the positive real axis using
the following notation.
\begin{equation*}
R_0^\pm(\lambda) := R_0(\lambda\pm i0) \qquad
R_V^\pm(\lambda) := R_V(\lambda\pm i0)
\end{equation*}

Every perturbation $V \in \Kato_0$ satisfies the local Kato condition 
\begin{equation} \label{eq:localKato}
\lim_{\delta\to 0} \sup_{y\in\R^3} \int_{|x-y| < \delta} \frac{|V(x)|}{|x-
y|}\,dx
 = 0.
\end{equation}
If $V$ has compact support then~\eqref{eq:localKato} is even equivalent
to membership in $\Kato_0$.
This degree of local regularity is sufficient to to conlude that
$H = -\Delta + V$ is essentially self-adjoint with spectrum bounded
below by $-M$ for some $M < \infty$~\cite{Si82}.  The Stone formula for the
absolutely continuous spectral measure of $H$ then dictates that
\begin{equation} \label{eq:Stone}
e^{-itH}f = \sum_j e^{-i\lambda_j t} P_{\lambda_j}(H)f
+ \frac{1}{2\pi i}\int_0^\infty e^{-it\lambda}
[R_V^+(\lambda)-R_V^-(\lambda)]f \,d\lambda.
\end{equation}
Because zero is assumed to be a regular point of the spectrum,
the summation contains only a finite set of eigenvalues
$\lambda_j$ which are removed by the projection $I - P_{pp}(H)$.  
Once the initial sum is
forced to vanish, dispersive estimates will succeed or fail based on the 
behavior
of the integral term.

It is customary to view the right-hand integral as a contour integral in the 
complex plane along a path that circles $\R^+$ clockwise.
Making a change of variables $\lambda \mapsto \lambda^2$ opens up the
contour to the entire real axis, with the understanding that
\begin{equation*}
R_V^+(\lambda^2) = \lim_{\eps\downarrow 0}R_V((\lambda + i\eps)^2)
 = \lim_{\eps\downarrow 0}R_V(\lambda^2 + i\,{\rm sign}\,(\lambda)\eps)
\end{equation*}
for all $\lambda \in \R$.  Written out this way the integral term in 
\eqref{eq:Stone} simplifies to
\begin{align*}
\frac{1}{\pi i}&\int_{-\infty}^\infty e^{-it\lambda^2} \lambda
R_V^+(\lambda^2)f \,d\lambda \\
&=\ \frac{1}{\pi i}\int_{-\infty}^\infty e^{-it\lambda^2} \lambda
\Res(\lambda^2)(I + V\Res(\lambda^2))^{-1}f\, d\lambda.
\end{align*}
A formal integration by parts leads to the expression
\begin{equation*}
\frac{1}{2\pi t}\int_{-\infty}^\infty e^{-it\lambda^2}
\big(I + \Res(\lambda^2)V\big)^{-1}\frac{d}{d\lambda}\Big[\Res(\lambda^2)
\Big] \big(I+V\Res(\lambda^2)\big)^{-1}f\, d\lambda.
\end{equation*}
If we adopt the shorthand notation $\hat{T}^\pm(\lambda) := 
VR_0^\pm(\lambda^2)$, 
%and also take $\hat{T}^-(\lambda) := I + VR_0^-(\lambda^2)$,
Theorem~\ref{thm:dispersive} should follow from
the estimate
\begin{multline} \label{eq:reduction}
\Big|\int_{-\infty}^\infty e^{-it\lambda^2} 
\Big\la \frac{d}{d\lambda}\Big[\Res(\lambda^2)\Big] 
(I+\hat{T}^+(\lambda))^{-1} f,\, (I+\hat{T}^-(\lambda))^{-1} g\Big\ra\, d\lambda 
\Big| \\
\les |t|^{-\frac12}
\norm[f][1] \norm[g][1].
\end{multline}
In fact Theorem~\ref{thm:dispersive} only requires~\eqref{eq:reduction}
to hold for all $f,g \in {\rm ran} (I-P_{pp}(H)) \subset L^1$.  The extra
restriction on the class of test functions is crucial when $H$ possesses
an eigenvalue at zero %(see~\cite{Go08p2} for dispersive estimates in this case)
but it is not needed here.

The choice of notation is designed to indicate that $\hat{T}^\pm(\lambda)$
exists as the Fourier transform of an important family of operators
$T^\pm(\rho)$, with $\rho$ serving as the variable dual to $\lambda$.
In particular, note that
the integral in \eqref{eq:reduction} can be estimated with the help of 
Plancherel's
identity.  Within the integrand there is a  factor of
$e^{-it\lambda^2}\frac{d}{d\lambda}[\Res(\lambda^2)]$. This family of
integral operators has an explicit kernel representation
\begin{equation*}
K(\lambda, x,y) = (-4\pi i)^{-1}e^{-it\lambda^2 + i \lambda |x-y|}.
\end{equation*}
Its inverse Fourier transform is also a family of convolution operators,
with kernel
\begin{equation*}
 \check{K}(\rho, x, y) 
  = (64\pi^3 i t)^{-\frac12} e^{-i\frac{(\rho + |x-y|)^2}{4t}}.
\end{equation*}
The expression for $\check{K}$ is bounded by $|t|^{-\frac12}$ for every value
of $\rho$, $x$, and $y$.
Therefore it suffices to prove that the inverse Fourier transforms of 
$(I + \hat{T}^+(\lambda))^{-1} f$ and 
$(I+\hat{T}^-(\lambda))^{-1}g$ both satisfy a corresponding
$L^1$ estimate. Our goal is then to prove the following mapping properties for 
$T^\pm(\rho)$.

\begin{theorem} \label{thm:goal}
Let $V \in \Kato_0$ be a scalar potential in $\R^3$ and suppose that 
$H = -\Delta + V$ has no resonances or eigenvalues along the interval
$[0,\infty)$.  Then
\begin{align*}
\norm[T^\pm(\rho)f][L^1_\rho L^1_x] \le 
  \bigg(\frac{\norm[V][\Kato]}{4\pi}\bigg)
  \norm[f][1]  \\
\text{and}\quad 
 \bignorm[\big((I+\hat{T}^\pm)^{-1}- I\big)^\vee (\rho)f
  ][L^1_\rho L^1_x] \le \ C \norm[f][1]
\end{align*}
for all $f \in L^1(\R^3)$.
\end{theorem}

The same approach is taken in~\cite{Go06b} where the potential is instead
assumed to belong to $L^{\frac32-\eps}(\R^3) \cap
L^{\frac32 + \eps}(\R^3)$.  The scale-invariant $L^p$ space for
Schr\"odinger potentials is $L^{\frac32}(\R^3)$ exactly, so the intersection
condition demands more decay at infinity and better local regularity than 
scaling arguments alone would suggest.  
We present Theorem~\ref{thm:goal} as an application of the operator-valued
Wiener Theorem in Section~\ref{sec:Wiener}.  This method has two principal
advantages over the previous work.  First, the result is sharper: The space
$\Kato_0$ includes all other known admissible classes of potentials
and the global Kato norm remains invariant under scaling transformations 
$V(x) \mapsto r^2V(rx)$.  Second, the proof is considerably cleaner as many of
the delicate $L^p$ resolvent bounds are replaced with a crude but effective
limiting argument.

For small potentials with $\norm[V][\Kato] < 4\pi$, the constant in the second
ineqality can be bounded by $(1 - \norm[V][\Kato]/4\pi)^{-1}$.  Dispersive estimates
for time-dependent potentials below this size threshold are 
proved in~\cite{RoSc04}.  When the potential is large the constant depends
more heavily on spectral properties (e.g.\ the avoidance of resonances)
and is not directly tied to the size of $V$.

\section{An Operator-Valued Wiener Theorem} \label{sec:Wiener}
Given a Banach space $X$, let $\W_X$ be the space of bounded linear
maps $T: X \to L^1(\R; X)$.  In the event that $X$ is not separable, the
target space $L^1(\R;X)$ is defined to be the norm-closure of simple functions.
%It is convenient, if somewhat inaccurate,
%to regard elements of $W_X$ as

$\W_X$ contains, for example, all measurable operator-valued
functions $T: \R \to \B(X)$
%(in fact $T(\rho)$ may be unbounded)
for which $\|T(\rho)\|_{\B(X)}$ is a.e.\ finite and the integral
\begin{equation}
\int_\R \norm[T(\rho)][\B(X)]\,d\rho
%\norm[T][\W_X] :=\sup_{\norm[f][X] = 1} \int_\R \norm[T(\rho)f][X]\,d\rho
\end{equation}
is also finite. However, not all elements of $\W_X$ have this form, as indicated in \cite{Be09p}.
 
$\W_X$ is an algebra under the formal convolution-composition product
\begin{equation} \label{eq:multiplication}
S * T(\rho)f = \int_\R S(\rho-\sigma)T(\sigma)f\, d\sigma
\end{equation}
whose boundedness follows from standard $L^1$ arguments including
approximating $T$ by an operator for which $Tf$ is a simple function.

There is a 
well-defined Fourier transform for elements of $\W_X$, namely
\begin{equation*}
\hat{T}(\lambda)f = \int_\R e^{-i\lambda\rho} T(\rho)f\, d\rho.
\end{equation*}
With no further assumptions on $T$ one can determine that $\hat{T}$ is
a strongly continuous family of operators with 
$\sup_{\lambda} \norm[\hat{T}(\lambda)][\B(X)] \le \norm[T][\W_X]$.
In addition $\hat{T}(\lambda)$ converges to zero in the
strong operator topology as $|\lambda| \to \infty$, by the Riemann-Lebesgue lemma.

There is no naturally occuring multiplicative identity in $\W_X$.  Let
$\overline{\W}_X$ denote the extension of $\W_X$ to include complex multiples
of an identity element, with the norm $\norm[z{\mathbf 1} + T][\overline{\W}_X]
= |z| + \norm[T][\W_X]$.  If necessary one can write ${\mathbf 1}(\rho) = 
\delta_0(\rho)I$, where $I$ is the identity in $\B(X)$, and its Fourier
transform is the constant function $\hat{\mathbf 1}(\lambda) = I$.
We note (as an aside) that $\W_X$ and
$\overline{\W}_X$ both embed isometrically as subalgebras of $\B(L^1(\R; X))$. 

As usual the Fourier transform intertwines convolution products in
$\overline{\W}_X$ with pointwise multiplication in 
$C^{\rm strong}_0(\R; \B(X))$.  For any pair of elements 
$z_1{\mathbf 1} + S$, $z_2{\mathbf 1} + T\in  \overline{\W}_X$ we have
\begin{equation}
\big[(z_1{\mathbf 1}+ S) * (z_2{\mathbf 1}+T)]\hat{\phantom{i}}(\lambda)
 = (z_1I + \hat{S}(\lambda))(z_2I + \hat{T}(\lambda)).
\end{equation}
Based on the product formula, any invertible element of $\overline{\W}_X$
must possess a Fourier transform that is invertible at every $\lambda \in \R$,
with uniformly bounded inverses.  An ideal Wiener theorem might show that these
conditions are sufficient for invertibility in $\overline{\W}_X$.  The
version that we prove here includes modest assumptions about the continuity
and locality of $T$.

\begin{theorem} \label{thm:Wiener}
Suppose $T$ is an element of $\W_X$ satisfying the
properties
\renewcommand{\theenumi}{(C\arabic{enumi})}
\renewcommand{\labelenumi}{\theenumi}
\begin{enumerate}
\item \label{translation} $\lim\limits_{\delta \to 0} 
\norm[T(\rho) - T(\rho-\delta)][\W_X] = 0$.
\item \label{locality}  $\lim\limits_{R \to \infty}
\norm[T\chi_{|\rho| \ge R}][\W_X] = 0$.
\end{enumerate}
If $I + \hat{T}(\lambda)$ is an invertible element of $\B(X)$ for every
$\lambda \in \R$, then ${\mathbf 1} + T$ possesses an inverse in 
$\overline{\W}_X$ of the form ${\mathbf 1} + S$.

In fact it is only necessary for some finite power $T^N \in \W_X$ (using the
definition of products in $\W_X$ given by~\eqref{eq:multiplication}) to satisfy
the translation-continuity condition~{\rm \ref{translation}} rather than $T$ 
itself.
\end{theorem}

\begin{proof}
It suffices to show that $(I + \hat{T}(\lambda))^{-1}$ is the Fourier transform
of an element ${\mathbf 1} + S \in \overline{\W}_X$. 
Let $\eta: \R \to \R$ be a standard
cutoff function.  For any real number $L$ one can express 
$(1 - \eta(\lambda/L))\hat{T}$ as the Fourier transform of
\begin{equation*}
S_L(\rho) = \big(T - L\check{\eta}(L\,\cdot\,) * T\big)(\rho) = 
\int_\R L \check{\eta}(L\sigma) [T(\rho) - T(\rho-\sigma)]\,d\sigma
\end{equation*}
Thanks to condition~\ref{translation}, the $\W_X$ norm of the
right-hand integral vanishes as $L \to \infty$.  This makes it possible
to construct an inverse Fourier transform for
\begin{equation*}
(1 - \eta(\lambda/2L))\big(I + \hat{T}(\lambda)\big)^{-1} 
= (1 - \eta(\lambda/2L)) 
 \sum_{k=0}^\infty (-1)^k \Big(\big(1 - \eta(\lambda/L))\hat{T}(\lambda)\Big)^k
\end{equation*}
via a convergent power series expansion in $S_L$ provided $L \ge L_1$. 
The $k=0$ term is exactly
${\mathbf 1} - 2L\check{\eta}(2L\rho)I$, and every subsequent term belongs
to $\W_X$.

If only $T^N$ satisfies~\ref{translation} then one constructs an
inverse Fourier transform for $(1-\eta(\lambda/2L))
 (I \pm \hat{T}^N(\lambda))^{-1}$ via the above process and observes that
\begin{equation*}
(1 - \eta(\lambda/2L))\big(I + \hat{T}(\lambda)\big)^{-1}
  = (1 - \eta(\lambda/2L))\big(I + (-\hat{T}(\lambda))^N \big)^{-1}
\sum_{k=0}^{N-1}(-1)^k \hat{T}^k(\lambda).
\end{equation*}

A similar approach works for finding a local inverse in the neighborhood of
any $\lambda_0 \in \R$.  For simplicity, consider the representative case 
$\lambda_0 = 0$, and let $A_0 = I + \hat{T}(0) \in \B(X)$.    
One can write $\eta(\lambda/L)(I + \hat{T}(\lambda) - A_0)$ 
as the Fourier transform of
\begin{align*}
S_L(\rho) &= L\check{\eta}(L\,\cdot\,) * T(\rho) - L\check{\eta}(L\rho)(I-A_0)
\\
&= \int_\R L\big(\check{\eta}(L(\rho-\sigma))- \check{\eta}(L\rho)\big)
 T(\sigma)\,d\sigma.
\end{align*}
Here we have used the fact that $A_0 = I + \int_\R T(\rho)\,d\rho$.

By the mean value theorem, $\int_\R L|\check{\eta}(L(\rho-\sigma))- 
 \check{\eta}(L\rho)|\,d\rho \les \min(L|\sigma|, 1)$. Recall that for
any fixed unit vector $f \in X$, the function $T(\sigma)f \in L^1(\R;X)$
has norm bounded by $\norm[T][\W_X]$.  If assumption~\ref{locality} holds then
there is $R_\eps$ so that $\norm[\chi_{|\sigma| > R_\eps} T(\sigma)f][] < 
\eps\norm[f][]$ uniformly in the choice of $f$.  It follows that
\begin{equation*}
\norm[S_L(\rho)f][L^1(\R;X)] \les (LR_\eps\norm[T][\W_X] + \eps)\norm[f][X]
\end{equation*}
and consequently that $\lim\limits_{L \to 0} \norm[S_L][\W_X] = 0$.

For any smooth function $\phi$ supported in $[-\frac{L}2,\frac{L}2]$, there
exists a local Neumann series
\begin{align*}
\phi(\lambda)(I + \hat{T}(\lambda))^{-1} 
&= \phi(\lambda)\big(A_0 + \eta(\lambda/L)(1 + \hat{T}(\lambda) - A_0)\big)^{-1}
\\
&= \phi(\lambda) A_0^{-1}\big(I + \hat{S_L}(\lambda)A_0^{-1}\big)^{-1}
 = \phi(\lambda) A_0^{-1} \sum_{k=0}^\infty (-1)^k 
      \big(\hat{S_L}(\lambda)A_0^{-1}\big)^k.
\end{align*}
So long as $L$ is chosen small enough that $\norm[S_L][\W_X]
\norm[A_0^{-1}][\B(X)] < 1$, the inverse Fourier transform of this series is 
convergent in the $\W_X$ norm.

Over the compact interval $\lambda_0 \in [-2L_1, 2L_1]$, there is a nonzero
lower bound on the length $L$ required for convergence of the above
power series.
Therefore it is possible to choose a partition of unity on $[-2L_1, 2L_1]$
so that the support of each cutoff $\phi_j$ has diameter small enough
that $\phi_j(\lambda)(I + \hat{T}(\lambda))^{-1}$ is the Fourier transform 
of an element of $\W_X$.  Finally, the inverse Fourier transform of
$1 - \eta(\lambda/2L)(I + \hat{T}(\lambda))^{-1}$ belongs to the affine
space ${\mathbf 1} + \W_X \subset \overline{\W}_X$, completing the
construction of $({\mathbf 1} + T)^{-1}$.
\end{proof}

\section{Proof of Theorem~\ref{thm:goal}}
We give the proof for $T^-(\rho)$.  Since $R_0^-(\lambda^2) = \Res((-
\lambda)^2)$
there is no difference between $T^+(\rho)$ and $T^-(\rho)$ except for a
reflection along the $\rho$-axis.

The first statement in Theorem~\ref{thm:goal} is a direct calculation.
Using the fact that in three dimensions $R_0^-(\lambda^2)$ is a convolution
against the kernel $(4\pi|x|)^{-1}e^{-i\lambda|x|}$ we can compute that
\begin{equation} \label{eq:Tminus}
T^-(\rho)f(x) = (4\pi \rho)^{-1} V(x) \int_{|x-y|= \rho} f(y)\,dy.
\end{equation}

Suppose $f$ is a bounded compactly-supported function in $\R^3$.
Then its convolution with the surface measure of a sphere is still bounded
and compactly supported.  Meanwhile $V \in \Kato_0$ is locally integrable
so $T^-(\rho)f$ belongs to $L^1(\R^3)$ for all $\rho \not= 0$.  Finally,
\begin{align*}
\int_{\R^3}\int_\R |T^-(\rho)f(x)|\, dx\,d\rho &\le
\frac{1}{4\pi} \int_{\R^3} \int_{\R^3} \frac{|V(x)|}{|x-y|}|f(y)|\,dy\,dx \\
&\le \frac{1}{4\pi} \norm[V][\Kato] \norm[f][1].
\end{align*}

The inequality extends by continuity to all $f \in L^1(\R^3)$, which shows
that $\norm[T^-][\W_{L^1}] \le \norm[V][\Kato]/4\pi$.  We remark once again
that the individual operators $T^-(\rho)$ need not belong to $\B(L^1(\R^3))$
for the integral inequality to hold.

The second inequality in Theorem~\ref{thm:goal} is the stronger conclusion.
It will be an immediate consequence of Theorem~\ref{thm:Wiener}, 
for $X = L^1(\R^3)$, once we
verify that $T^-$ satisfies condition~\ref{locality} and $(T^-)^N$
satisfies~\ref{translation} for some finite $N$.  Pointwise invertibility
of $I + \hat{T}^-(\lambda)$ comes from the Fredholm alternative and the
assumption that $H$ has no eigenvalues or resonances along its
continuous spectrum.  

To summarize the invertibility argument:
$\hat{T}^-(\lambda) = V R_0^-(\lambda^2)$ is a compact operator
on $L^1(\R^3)$ for each $\lambda$, improving regularity by two derivatives,
if $V$ belongs to $C^\infty_c$.  Compactness is preserved by
taking norm limits to $V \in \Kato_0$. By the Fredholm alternative,
$I + \hat{T}^-(\lambda_0)$ fails to be invertible only if there exists a
function $\vp \in L^{1}$ such that $(I + \hat{T}^-(\lambda_0))\vp = 0$. 
Then $\psi = R_0^-(\lambda_0^2) \vp$ would be a solution to
$\psi + R_0^-(\lambda_0^2) V \psi = 0$ satisfying 
$(1+|x|)^{-s}\psi \in L^2(\R^3)$ for all $s > \frac12$ because the free
resolvent is a bounded map from $L^1$ into these weighted function spaces.

%Furthermore, $\hat{T}^-(\lambda_0)$ being compact or $I + \hat{T}^-(\lambda_0)$ being invertible on $L^1(\R^3)$ is
%equivalent to $\hat{T}^-(\lambda_0)^* = R_0^+(\lambda^2) V$ being compact, respectively $I +
%\hat{T}^-(\lambda_0)^*$ being invertible on $L^{\infty}(\R^3)$. By the Fredholm alternative, then, the lack of invertibility
%implies the existence of an $L^{\infty}$ solution $\psi$ to the equation $\psi + R_0^-(\lambda_0^2) V \psi = 0$.

Either way, $H$ would have an eigenvalue or resonance at $\lambda_0^2$
according to whether or not $\psi \in L^2(\R^3)$.

The locality condition~\ref{locality} is rather straightforward as well.
Suppose $V$ is a bounded function with compact support in a set of
diameter $D$.  From~\eqref{eq:Tminus} it follows that for $R > 2D$
\begin{align*}
\int_{\R^3}\int_{|\rho|>R} |T^-(\rho)f(x)|\, dx\,d\rho &\le
\frac{1}{4\pi} \iint_{|x-y|>R} \frac{|V(x)|}{|x-y|}|f(y)|\,dy\,dx \\
&\les R^{-1} \norm[V][1] \norm[f][1]
\end{align*}
Then $\norm[T^-\chi_{|\rho| > R}][\W_{L^1}] \to 0$ as $R \to \infty$,
and this property is preserved under a limiting sequence in $V \in \Kato$.

\begin{remark}
Condition~\ref{locality} is actually satisfied for the larger class of 
potentials that possess the "distal Kato property"
\begin{equation*} 
\lim_{R\to \infty} \sup_{y \in \R^3} 
  \int_{|x-y| > R} \frac{|V(x)|}{|x-y|}\,dx = 0
\end{equation*}
Unlike the assumption $V \in \Kato_0$, the distal Kato property does
not guarantee that operators $VR_0^-(\lambda)$ act compactly
on $L^1(\R^3)$.  Most technical elements of the dispersive
estimate are not affected except for the question of whether
$I + VR_0^-(\lambda)$ is invertible pointwise in $\lambda$.
Normally one uses the Fredholm Alternative to derive it from the
absence of embedded eigenvalues and resonances.  With this argument
unavailable one needs to strengthen the spectral assumptions on
$H$ accordingly.
\end{remark}

The remaining task is to verify that some power of $T^-$ in $\W_{L^1}$, $(T^-)^N$, satisfies the
translation-continuity hypothesis~\ref{translation}.  We proceed using the
value $N=4$. Choose a bounded compactly supported approximation $V_\eps$
with $\norm[V-V_\eps][\Kato] < \eps$ and $\norm[V_\eps][\Kato] \le 
\norm[V][\Kato]$.  Let $T_\eps^-$ denote the corresponding element of 
$\W_{L^1}$.

As before the compact support makes it possible to choose $R$ so that 
$\norm[T_\eps^- \chi_{|\rho| > R}][\W_{L^1}] $ is smaller than
$ \eps\norm[V][]^{-3}$.  Then by convolution of support (in $\rho$) one also 
has $\norm[(T_\eps^-)^4\chi_{|\rho| > 4R}][\W_{L^1}] < 4\eps$, with the
result that
\begin{equation*}
\bignorm[(T_\eps^-)^4(\rho)-(T_\eps^-)^4(\rho-\delta)][\W_{L^1}]
\le \bignorm[((T_\eps^-)^4(\rho)-(T_\eps^-)^4(\rho-\delta))\chi_{|\rho| < 4R+1}
 ][\W_{L^1}] + 8\eps
\end{equation*}
for any $\delta < 1$.  Over the finite interval $|\rho| < 4R + 1$ it will
suffice to show that $(T_\eps^-)^4(\rho)$ is in fact a continuously
differentiable function taking values in $\B(L^1(\R^3))$.

The kernel of the free resolvent has the pointwise absolute value
$|R_0^-(\lambda^2){\scr (x,y)}| = \frac{1}{4\pi|x-y|}$
which shows it to be a bounded operator from $L^1(\R^3)$ to 
$L^2({\rm supp}(V_\eps))$.  Meanwhile there is also a known
family of weighted $L^2$ estimates
\begin{equation*}
\norm[\japanese[x]^{-\alpha}R_0^-(\lambda^2)\japanese[x]^{-\alpha}f][2]
  \le C_\alpha (1 + |\lambda|)^{-1}\norm[f][2]
\end{equation*}
for any exponent $\alpha > \frac12$ (cf.~\cite[Theorem~5.1]{AgHo76}).
Taken in combination these imply that
\begin{equation*}
\bignorm[\big((T_\eps^-)^4\big)\hat{\phantom{i}}(\lambda) f][1] =
\bignorm[\big(V_\eps R_0^-(\lambda^2)\big)^4 f][1]  \les 
    (1+|\lambda|)^{-3}\norm[f][1]
\end{equation*}
with the constant determined by the maximum size of $V_\eps$ and the diameter
of its support.  With this much decay present as $\lambda \to \infty$ it follows
from Fourier inversion that $(T_\eps^-)^4(\rho)$ has a derivative bounded in 
size
by the same constant.  Therefore
\begin{equation*}
\bignorm[((T_\eps^-)^4(\rho)-(T_\eps^-)^4(\rho-\delta))\chi_{|\rho| < 4R+1}
 ][\W_{L^1}] \le C(V_\eps) (4R+1) \delta
\end{equation*}
and $\delta$ can be chosen sufficiently small to keep this quantity less
than $\eps$.  The norm difference between $(T_\eps^-)^4$ and its
translation is no greater than $9\eps$.

The triangle inequality permits a step from $T_\eps^-$ back to $T^-$
by giving the bound
\begin{align*}
\bignorm[(T^-)^4(\rho) - (T^-)^4(\rho-\delta)][\W_{L^1}]
  &\le 2\norm[(T^-)^4 - (T_\eps^-)^4][\W_{L^1}] + 9\eps \\
  &\le (C\norm[V][\Kato]^3 + 9)\eps
\end{align*}
for all $\delta$ sufficiently small.  Taking $\eps$ to zero verifies that
condition~\ref{translation} is satisfied.

%\nocite{BuPlStTa03}
%\nocite{Du07}
%\nocite{Je80}
%\nocite{Je84}
%\nocite{Ka65}
%\nocite{Ra78}
%\nocite{ReSi4}

\bibliographystyle{abbrv}
\bibliography{MasterList}

\begin{thebibliography}{10}

\bibitem{AgHo76}
S.~Agmon and L.~H\"ormander.
\newblock Asymptotic properties of solutions of differential equations with
  simple characteristics.
\newblock {\em J.\ Anal.\ Math.}, 30(1):1--38, 1976.

\bibitem{Be09p}
M.~Beceanu.
\newblock New estimates for a time-dependent {S}chr\"odinger equation.
\newblock arXiv:0909.4029, 2009.

\bibitem{BeGo10p}
M.~Beceanu and M.~Goldberg.
\newblock Scaling-invariant {S}trichartz estimates for the wave equation.
\newblock In preparation.

\bibitem{BuPlStTa03}
N.~Burq, F.~Planchon, J.~Stalker, and A.~S. Tahvildar-Zadeh.
\newblock Strichartz estimates for the wave and {S}chr\"odinger equations with
  the inverse--square potential.
\newblock {\em J.\ Funct.\ Anal.}, 203(2):519--549, 2003.

\bibitem{Go06b}
M.~Goldberg.
\newblock Dispersive bounds for the three-dimensional {S}chr\"odinger equation
  with almost critical potentials.
\newblock {\em Geom.\ and Funct.\ Anal.}, 16(3):517--536, 2006.

\bibitem{GoSc04a}
M.~Goldberg and W.~Schlag.
\newblock Dispersive estimates for the {S}chr\"odinger operator in dimensions
  one and three.
\newblock {\em Comm.\ Math.\ Phys.}, 251(1):157--178, 2004.

\bibitem{GoVeVi06}
M.~Goldberg, L.~Vega, and N.~Visciglia.
\newblock Counterexamples of {S}trichartz inequalities for {S}chr\"odinger
  equations with repulsive potentials.
\newblock {\em Intl.\ Math.\ Res.\ Not.}, 2006:16pp., 2006.
\newblock Article ID 13927.

\bibitem{JoSoSo91}
J.-L. Journ\'e, A.~Soffer, and C.~Sogge.
\newblock Decay estimates for {S}chr\"odinger operators.
\newblock {\em Comm.\ Pure Appl.\ Math.}, 44(5):573--604, 1991.

\bibitem{RoSc04}
I.~Rodnianski and W.~Schlag.
\newblock Time decay for solutions of {S}chr\"odinger equations with rough and
  time-dependent potentials.
\newblock {\em Invent.\ Math.}, 155(3):451--513, 2004.

\bibitem{Si82}
B.~Simon.
\newblock Schr\"odinger semigroups.
\newblock {\em Bull.\ Amer.\ Math.\ Soc.}, 7(3):447--526, 1982.

\bibitem{Ya95}
K.~Yajima.
\newblock {The $W^{k,p}$-continuity of wave operators for {S}chr\"odinger
  operators}.
\newblock {\em J. Math.\ Soc.\ Japan}, 47(3):551--581, 1995.

\end{thebibliography}

\end{document}